# A METHOD OF SOLVING CERTAIN NONLINEAR DIOPHANTINE EQUATIONS

Florentin Smarandache, Ph D
Professor of Mathematics
Department of Math & Sciences
University of New Mexico
200 College Road
Gallup, NM 87301, USA
E-mail: smarand@unm.edu

**Abstract**.
In this paper we propose a method of solving a Nonlinear Diophantine Equation by converting it into a System of Diophantine Linear Equations.

**Introduction**.
Let's consider a polynomial with integer coefficients, of degree $m$

$$P(X_1,...,X_n) = \sum_{\substack{0 \le i_1+...+i_n \le m \\ 0 \le i_j \le m, j=1,n}} a_{i_1...i_n} X_1^{i_1}...X_n^{i_n}$$

which can be factored out in linear factors (that can eventually be established through the undetermined coefficients method):

$$P(X_1,...,X_n) = \left(A_1^{(1)}X_1 + ... + A_n^{(1)}X_n + A_{n+1}^{(1)}\right) \cdots \left(A_1^{(m)}X_1 + ... + A_n^{(m)}X_n + A_{n+1}^{(m)}\right) + B$$

with all $A_j^{(k)}$, $B$ in $\mathbb{Q}$, but which by bringing to the same common denominator and by eliminating it from the equation $P(X_1,...,X_n) = 0$ they can be considered integers.. Thus the equation transforms in the following system:

$$\begin{cases} A_1^{(1)}X_1 + ... + A_n^{(1)}X_n + A_{n+1}^{(1)} = D_1 \\ ............................................. \\ A_1^{(m)}X_1 + ... + A_n^{(m)}X_n + A_{n+1}^{(m)} = D_m \end{cases}$$

where $D_1,...,D_m$ are the divisors for $B$ and $D_1 \cdots D_m = B$.

We separately solve each linear Diophantine equation and then we intersect the equations' solutions.

**Example 1.**
Solve in integer numbers the equation:
$$-2x^3 + 5x^2y + 4xy^2 - 3y^3 - 3 = 0.$$
We'll write the equation in another format
$$(x+y)(2x-y)(-x+3y) = 3.$$
Let $m$, $n$ and $p$ be the divisors of $3$, $m \cdot n \cdot p = 3$. Thus

$$\begin{cases} x + y = m \\ 2x - y = n \\ -x + 3y = p \end{cases}$$

For this system to be compatible it is necessary that

$$\begin{pmatrix} 1 & 1 & m \\ 2 & -1 & n \\ -1 & 3 & p \end{pmatrix} = 0,$$

or

$$5m - 4n - 3p = 0 \quad (1)$$

In this case

$$x = \frac{m+n}{3} \text{ and } y = \frac{2m-n}{3} \quad (2)$$

Because $m, n, p \in \mathbb{Z}$, from (1) it results – by solving in integer numbers – that:

$$\begin{cases} m = 3k_1 - k_2 \\ n = k_2 \\ p = 5k_1 - 3k_2 \end{cases} \quad k_1, k_2 \in \mathbb{Z}$$

which substituted in (2) will give us $x = k_1$ and $y = 2k_1 - k_2$. But $k_2 \in D(3) = \{\pm 1, \pm 3\}$; thus the only solution is obtained for $k_2 = 1$, $k_1 = 0$ from where $x = 0$ and $y = -1$.

**Example 2**.

Analogously, it can be shown that, for example the equation:
$$-2x^3 + 5x^2y + 4xy^2 - 3y^3 = 6$$
does not have solutions in integer numbers.

**REFERENCES**

[1] Marius Giurgiu, Cornel Moroti, Florică Puican, Stefan Smărăndoiu – Teme şi teste de Matematică pentru clasele IV-VIII - Ed. Matex, Rm. Vîlcea, Nr. 3/1991

[2] Ion Nanu, Lucian Tuţescu – "Ecuaţii Nestandard", Ed. Apollo şi Ed. Oltenia, Craiova, 1994.